\documentclass[11pt]{article}
\usepackage{amsmath}
\usepackage{amsthm}
\usepackage{amsfonts}

\title{Discrete spectral triples and their symmetries}
\author{
Mario Paschke\thanks{e-mail:paschke@dipmza.physik.uni-mainz.de} \\
Andrzej Sitarz\thanks{e-mail:sitarz@higgs.physik.uni-mainz.de} \\ \ \\
Universit\"at,  55130 Mainz, Germany }
\begin{document}

\begin{titlepage}

\maketitle

{\bf Abstract:}
We classify $0$-dimensional spectral triples over complex and
real algebras and provide some general statements about their
differential structure. We investigate also whether such spectral
triples admit a symmetry arising from the Hopf algebra structure of the
finite algebra. We discuss examples of commutative algebras and group
algebras. 

\vskip 2cm

\begin{flushright}
MZ-TH/96-40 \\
q-alg/yymmddd \\
\end{flushright}

\end{titlepage}

\parindent 5pt
\def\id{\hbox{id\ }}
\def\ker{\hbox{ker\ }}
\def\CA{{\cal A}}
\def\CH{\hbox{$\cal H$}}
\def\CJ{{\cal J}}
\def\C{\mathbb C}
\def\R{\mathbb R}
\def\H{\mathbb H}
\def\Z{\mathbb Z}
\def\bo{\hbox{\bf 1}}
\def\zero{\hbox{\Large \ \ 0 \ \ }}
\def\Om#1{\Omega^{#1}(\CA)}
\def\tsa{\otimes_{\CA}}
\def\mtw#1#2#3#4{%
    \begin{array}{cc} #1 & #2 \\ #3 & #4 \end{array}}
\def\mth#1#2#3#4#5#6#7#8#9{%
    \begin{array}{ccc} #1 & #2 & #3 \\ #4 & #5 & #6 \\
                        #7 & #8 & #9  \end{array}}
\font\mathbbfont=msam10
\def\Box{\hbox{\mathbbfont \char"4}}
\newtheorem{lem}{Lemma}
\newtheorem{obs}{Observation}
\section{Introduction}

Recently Connes \cite{Connes} has proposed a definition of {\em
spectral triples}, which is supposed to extend the notion of the
Riemannian geometry to a noncommutative setup. This could be a natural
geometrical framework for theories of fundamental interactions and
gravity at a small distance scale and possibly it could open a way
towards quantum gravity. 

One of the promising applications to fundamental particle physics is an
attempt to explain the particle content and interaction pattern in the
Standard Model. In the spectral triple in question the finite algebra 
$M_3(\C) \oplus \H \oplus \C$ plays the significant role. Here we shall not 
discuss this particular physical application. \footnote{For a
detailed discussion of the Standard Model spectral triple see our
forthcoming paper \cite{PS2}.}.

In this paper we derive the general structure of $0$-dimensional spectral
triples and of the Dirac operator for complex and real semisimple
algebras. Moreover, assuming that such algebras are Hopf algebras we
discuss possible symmetry structures. As examples we present $\C^n$
commutative algebra with discrete group structure as well as group
algebras, particularly $\C{}S_3$ for the permutation group $S_3$.

\section{Complex Spectral Triples}

Let  $\CA$ be a finite dimensional semi-simple
algebra\footnote{$C^\ast$ algebras are semisimple, thus we consider
only them. Of course, nonsemisimple algebras (grassman variables, for
instance) may also appear in physical models. Whether one can
successfully construct a corresponding theory for them remains an
interesting problem.} over the field of complex numbers $\C$. Such an
algebra is isomorphic to $M_{n_1}(\C) \oplus \cdots \oplus 
M_{n_k}(\C)$ for some non-zero integers $n_1,\ldots n_k$ (see
\cite{EMS}, pp.29-36). 

A {\em spectral triple} for a finite algebra means 
a finite Hilbert space $\CH$, 
a faithful representation $\pi: \CA \to B(\CH)$ over the field $\C$, a
reality structure $J$, such that $J$ is an antilinear isometry of
$\CH$, $J^2=1$ and $\pi^0 = J \circ \pi^* \circ J$ is a representation
of the opposite algebra $\CA^0$, which commutes with $\pi$: 

$$  [\pi^0(a), \pi(b)]=0, \;\;\;\; \forall a,b \in \CA $$

We assume that there exists a grading $\gamma \in B(\CH)$, such that
$\gamma^* = \gamma$, $\gamma^2=1$, $J \gamma = \gamma J$, $ \gamma
\pi(a) = \pi(a) \gamma$ and $\gamma$ is equal to $\pi(c)$ for some
Hochschild cocycle $c \in Z_0(\CA,\CA \otimes \CA^0)$.

We begin with some simple observations. 

\begin{obs}
$\CH$ is a bimodule over $\CA$. If $\xi \in \CH$ then we have: 
$a\xi = \pi(a)\xi$ and $\xi a = \pi^0(a) \xi$.
\end{obs}

We can now use the structure of the finite algebra $\CA$ to learn more
about the Hilbert space of the spectral triple.

\begin{lem} Let $P_i$ be the element of $\CA$ with identity matrix in
the i-th entry and zeroes elsewhere. Let $H_{ij}$ be $P_i \CH P_j$.
Then $\CH$ is a direct sum of all $H_{ij}$:

$$ \CH = \bigoplus_{i,j} H_{ij} $$

and if for every $a \in \CA$ we denote $P_i a$ ($= a P_i$) as $a_i$ we
have:

$$ (a \xi)_{ij} = a_i \xi_{ij}, \;\;\;\; (\xi a)_{ij} = \xi_{ij} a_j $$

\end{lem}

The proof is a simple calculation and we skip it.

\begin{lem}

The Hilbert subspace $H_{ij}$ has the following form:

$$ H_{ij} = \C^{n_i} \otimes \C^{r_{ij}} \otimes \C^{n_j} $$

where $r_{ij}$ are some integers, $r_{ij} \geq 0$. We take $H_{ij}=0$
iff $r_{ij}=0$. An element $a$ acts on $H_{ij}$ from the left as $a_i
\otimes 1 \otimes 1$ and from the right as $1 \otimes 1 \otimes a_j^T$,
where ${}^T$ denotes the matrix transposition. \end{lem}

{}From the previous lemma we know that $H_{ij}$ is the representation
space for $M_{n_i}(\C)$ and for the opposite algebra of $M_{n_j}(\C)$.
Because these representations commute, the structure of $H_{ij}$ must
have the above form.

\begin{lem}
The grading operator $\gamma$ is $\pm 1$ when restricted to the 
subspace $H_{ij}$. 
\end{lem}

{\em Proof:} Since $\gamma$ commutes both with the elements of $\CA$
and $J$, and is self-adjoint, we conclude that it also commutes with
the right multiplication on $\CH$:

$$ a (\gamma \xi) b = \gamma (a \xi b) $$

Hence, $\gamma H_{ij} \subset H_{ij}$. When restricted to $H_{ij}$,
$\gamma$ is an isometry, which commutes with the left multiplication by
matrices from $M_{n_i}(\C)$ and right multiplication by matrices from
$M_{n_j}(\C)$. Therefore it must be of the form:

$$ \gamma_{ij} : H_{ij} \to H_{ij}, \;\; \gamma = 1_{n_i} \otimes
\Gamma_{ij} \otimes 1_{n_j}, $$

where $\Gamma_{ij}$ is an arbitrary self-adjoint matrix, such that
$\Gamma_{ij}^2=1$ acting on $\C^{r_{ij}}$. However, if we take into
account that $\gamma$ is an image of a Hochschild cocycle, i.e., there
exists an element $c$ of $\CA \otimes \CA^0$ such that $\pi(c)=\gamma$,
we see that $\Gamma_{ij}$ must be proportional to the identity matrix,
so it leaves only the possibility for $+1$ or $-1$.\Box

\begin{lem} The reality $J$ maps $H_{ij}$ onto $H_{ji}$, therefore
$r_{ij}$ must be equal to $r_{ji}$. \end{lem}

{\em Proof}. Let us take $\xi \in H_{ij}$. There exists $v \in \CH$
such that $\xi = P_i v P_j $. Next, calculate $J \xi$, writing all 
expressions carefully:

\begin{eqnarray*}
& J\xi = J P_i v P_j = J \pi(P_i) \left( J \pi(P_j) J \right) v = \\ &
= ( J P_i J) P_j (J v) = P_j (J v) P_i \in H_{ji} \end{eqnarray*}

Of course, since $J^2 = \id$ we must have $r_{ij}=r_{ji}$ so that the
dimensions of both Hilbert spaces are equal.\Box 

\begin{obs} If we denote by $\gamma_{ij}$ the sign of the grading on
$H_{ij}$ then $\gamma_{ij} = \gamma_{ji}$. \end{obs}

This is a trivial consequence of the commutation relation
$[J,\gamma]=0$.

We can now introduce $q_{ij} = r_{ij} \gamma_{ij}$, which, later on,
turns out to be the intersection form for our spectral triple.
So far, we know that $q_{ij}$ is a symmetric matrix with integer
entries.

\begin{lem} There exists a basis of $H_{ij}, H_{ji}$ of the form 
$ H_{ij} \ni v = v_i \otimes v_{ij} \otimes v_j$ such that  
\begin{equation}
Jv = v_j \otimes  v_{ji} \otimes v_i  \in H_{ji}. \label{J-prop}
\end{equation}
\end{lem}

{\em Proof:} We begin with the $i \not= j$ case. First of all, let us
fix an orthonormal basis of the spaces $H_{ij}$ and $H_{ji}$ of the
form $e_1 \otimes e_{2} \otimes e_3$ and define $\tilde{J}$ as an
antilinear isometry operator, which exchanges the first and the third
element of the tensor product:

$$ \tilde{J} ( e_1 \otimes e_{2} \otimes e_3 ) =
e_3 \otimes e_{2} \otimes e_1 $$

Then $\tilde{J}$ satisfies all the axioms for the spectral triple,
i.e., $\tilde{J}^2=1$ and $a^0 = \tilde{J} a^\dagger \tilde{J}$.

Next, observe that $ J \circ \tilde{J}$ is a linear invertible map,
such that: $$ J \circ \tilde{J} a = a J \circ \tilde{J} $$ and $$ J
\circ \tilde{J} a^0 = a^0 J \circ \tilde{J} $$ We prove only the first
of the two: \begin{eqnarray*} &  J \circ \tilde{J} a =  J \circ
\tilde{J} a \tilde{J}^2 = J (a^\dagger)^0 \tilde{J} =  \\ & = J^2 a J
\circ \tilde{J} =  a J \circ \tilde{J}. \end{eqnarray*}

Therefore $J \circ \tilde{J}$ commutes both with the left and the
right multiplication, so it has to be of the form: $1 \otimes \tilde{l}
\otimes 1$, where $\tilde{l}$ is a linear isometry  $l: \C^{r_{ij}} \to
\C^{r_{ji}}$. Thus $J = (1 \otimes \tilde{l}^{-1} \otimes 1) \circ
\tilde{J}$ and the only possible nontrivial part of $J$ comes from
$\tilde{l}$.

In the last step we have to consider the cases $j=i$ and $j\not=i$
separately. The latter one is simpler, we just chose any orthonormal
basis $e_s$ of $C^{r_{ij}}$ (same as in the definition of $\tilde{J}$)
and then $\tilde{l}^{-1} e_s$ constitutes  a basis of $C^{r_{ji}}$.
Take as $v_{ij}$ just the vectors of the chosen basis and $v_{ji}$
their image under $\tilde{l}^{-1}$. Then the property (\ref{J-prop})
follows immediately.

A more subtle situation occurs in the $i=j$ case.  Here, however, we
can use another simple argument. Take $v \in C^{r_{ii}}$. Then either
$\tilde{l}^{-1} v$ is proportional to $v$ or is linearly independent
from it. In the latter case it is sufficient to consider $v_{ii}$ to be
$v + \tilde{l}^{-1} v$ and $ i( v - \tilde{l}^{-1} v) $. In the first
case $\tilde{l}^{-1} v = e^{i\phi} v$ and then $e^{\frac{i\phi}{2}} v $
will do as $v_{ii}$. Since our space is finite dimensional and $l$ is
an isometry we can carry on with this procedure until we find a
complete basis of $H_{ii}$ satisfying our requirement.\Box

\subsection{The Dirac operator}

Up to this point we have constructed the basic geometry of a finite
spectral triple and in the next step we shall introduce the Dirac
operator, which gives the differential and metric structures.

The Dirac operator $D$ of a finite spectral triple is a self-adjoint
linear operator on $\CH$, which commutes with $J$: $DJ=JD$,
anticommutes with $\gamma$: $D\gamma=-\gamma D$ and satisfies the
following relation: \begin{equation} \left[ [ D, \pi(a)], \pi^0(b)
\right] = 0, \;\;\; \forall a,b \in \CA. \label{D-prop} \end{equation}

\begin{obs} Let $D_{ij,kl} = P_{ij} D P_{kl}$, where $P_{ij}$ is the
projection operator on $H_{ij}$: \mbox{($P_{ij} v = P_i v P_j$ )}. Then
the following is true:

1. $D_{ij,kl} : H_{kl} \to H_{ij}$ is a linear map, which is zero
unless $\gamma_{ij} \gamma_{kl} < 0$ (i.e., they are of opposite sign)

2. $D_{ij,kl}^\dagger = D_{kl,ij}$. This follows easily from the 
requirement $D=D^\dagger$.

3. $D_{ij,kl} = J D_{ji,lk} J $. This follows from the relation
$[D,J]=0$ and $J^2=\id$. If we choose the basis of the spaces $H_{ij}$
so that the operator $J$ has a preferred form $J v = \pm v$, the above
condition can be simplified to $D_{ij,kl} = D_{ji,lk}^*$, where $*$
means the complex conjugation of matrix elements of $D_{ji,lk}$ in the
preferred basis. \end{obs}

The strongest restriction for the Dirac operator follows from the
commutation property (\ref{D-prop}). 

\begin{lem}
The Dirac operator component $D_{ij,kl}$ vanishes unless $i=k$ or
$j=l$. In either of these cases $D$ must commute with the corresponding
action of the $\CA$ or $\CA^0$, respectively. \end{lem}

{\em Proof:} Let us rewrite the condition (\ref{D-prop}) for $v \in
H_{kl}$: $$ D_{ij,kl} (a_k v b_l) - a_i D_{ij,kl} (v b_l)
 - \left( D_{ij,kl} (a_k v) \right) b_j + a_i (D_{ij,kl} v) b_j = 0 $$
Choosing $b = P_j$ $j\not=l$ we get that $D_{ij,kl}=0$ unless $i=k$
and, for $i=k$, $D_{ij,il} a_i = a_i D_{ij,il}$ for every $a \in \CA$. 

Similarly, if we choose $a=P_i$, $i \not=k$ we get that 
$D_{ij,kl}=0$ unless $j=l$ and additionally $D_{ij,kj} (v b_j) = 
( D_{ij,kj} v) b_j$.

Of course, the case $i=k$ and $j=l$ is excluded by the requirement that
$D$ acts only between spaces of different grading.\Box

\subsection{Differential algebra}

The Dirac operator $D$ provides us with both the metric and the
differential properties on the finite space corresponding to the
algebra $\CA$. We shall now concentrate on some general properties of
the differential structure.

\begin{obs}
The representation $\pi$ of $\CA$ extends to the universal differential
algebra $\Omega_u\CA$, by: $$ \pi( a_0 da_1 d_2 \cdots da_n) = \pi(a_0)
[D, \pi(a_1)] \cdots [D,\pi(a_n)].$$ \end{obs}

The differential calculus is usually constructed as follows: first,
$\Omega^1(\CA)$ is the quotient of $\Omega_u^1\CA$ by the kernel of
$\pi$. Therefore the bimodule $\Om{1}$ is isomorphic to
$\pi(\Omega^1_u(\CA))$. Higher order forms are obtained by taking the
quotient of the universal forms of a given order by the ideal
$\Omega_u^n\CA \cap (\ker \pi  \cup d \ker \pi)$.

Let us state some general results.

\begin{lem} The first order differential calculus is inner, there
exists an one-form: \begin{equation} \xi = \sum_{i\not= j} P_i \, dP_j,
\end{equation} such that for every $a \in \CA$ one has: $da = [\xi,
a]$. The Dirac operator $D = \pi(\xi) + J \pi(\xi) J$. \end{lem}

{\em Proof:} Of course, from the definition we have $\pi(da) = [D,
\pi(a)]$, so one has to show:

$$ \sum_{i \not= j} \left[ \pi(P_i) \left[ D, \pi(P_j) \right], 
\pi(a) \right] = \left[ D, \pi(a) \right]. $$

Let us compute $\xi_{ij,kl} = P_{ij} \pi(\xi) P_{kl}$. For
any $h \in H_{kl}$ we have:
\begin{eqnarray*}
\left(\pi(\xi) h \right)_{ij} & = & 
\left( \sum_{s \not= t} P_s \left[ D, P_t \right] h \right)_{ij} \\
 & = & \sum_{s \not=t} (P_s D P_t h)_{ij} 
  =  (P_i D P_k h)_{ij} = \delta_{jl} D_{ij,kj} h.
\end{eqnarray*}
Note that:
$$ \xi = \sum_{i \not= j} P_i \, dP_j = - \sum_i P_i \, dP_i. $$

The relation $D= \pi(\xi) + J \pi(\xi) J$ follows immediately from the
properties of $J$ and $D$, additionally $[ J\xi J, a] =0$ for every
element  $a \in \CA$, so that the property $da = [D,a] =[\xi,a]$
holds.\Box

\begin{obs}
$\pi(\xi)$ is a selfadjoint operator, thus the first order calculus has
a natural involution map such that $d \circ \ast = - (\ast \circ\, d)$
if we define $\xi^\ast = \xi$. \end{obs}

The one-form $\xi$ has some more interesting properties, which make the
task of determining the differential structure easier. Before we
discuss the construction  of $\Om{2}$ let us make another useful
remark.

\begin{obs} 
The bimodule $\Om{1}$ has no center, if $a\omega=\omega a$ for every
element $a \in \CA$ then $\omega=0$. \end{obs}

{\em Proof:} First, let us observe that for every idempotent $e=e^2$,
we have $ e\, de\, e =0$. This follows from differentiating $0=e^2-e$
and multiplying the result by $e$  from the left and from the right. Therefore
$P_i dP_i P_i = P_i \xi P_i =0$.

Every element of $\Om{1}$ could be written as a finite sum
of the type $\omega = \sum_\alpha b_\alpha \xi c_\alpha$. Suppose that
it commutes with every $a \in \CA$, in particular with $P_i$. Then if
we multiply both sides by $P_i$, since $P_i^2=P_i$ we get $P_i \omega
P_i = P_i \omega$. However $P_i \omega P_i =0$ as $P_i$ commute with
elements of the algebra and $P_i \xi P_i=0$. So $P_i \omega=0$ for
every $i$, and therefore $\omega=0$.\Box

\begin{lem} 
$d\xi = \xi\xi + \sum P_i \xi \xi P_i $
\end{lem}
{\em Proof:}  
\begin{eqnarray*}
d \xi & = & - \sum_{i} dP_i dP_i \\
& = & - \sum_{i} [\xi, P_i] [\xi, P_i] 
 =   \xi \xi + \sum_i P_i \xi \xi P_i. \Box
\end{eqnarray*}
The above equality holds only within the differential algebra,
which means that all products of one-forms in the line above 
must be calculated using the product in $\Om{2}$. We have not
determined it yet (in a general situation it is a difficult
technical calculation), however, using the above relations we will be
able to say something about it.

\begin{lem} Let $\Xi \in \Om{1} \tsa \Om{1} \ni \Xi$ be
\begin{equation}
\Xi = \sum_i P_i \xi \tsa \xi P_i.
\end{equation}
 Let $\CJ_1 = \ker \pi$ (where $\pi$
on $\Om{1} \tsa \Om{1}$ is defined in the usual way:
$\pi(\omega_1 \tsa \omega_2) = \pi(\omega_1)\pi(\omega_2)$),
$\CJ_2$ be the subbimodule generated by commutators $[a, \Xi]$ for all 
$a \in \CA$ and finally let $\CJ_3$ be the subbimodule generated by all
elements of the form $\sum_i a_i (\xi \tsa \xi - \Xi) b_i$ for all
$a_i,b_i$ such that $\sum_i a_i \xi b_i = 0$. Then $\CJ = \CJ_1 \oplus
\CJ_2 \oplus \CJ_3$ is a subbimodule of $\Om{1} \tsa \Om{1}$ and
$\Om{2} = \Om{1}\tsa \Om{1} / \CJ $. \end{lem}

{\em Proof:} We know that every element of $\Om{1}$ can be written as a
finite sum of $\omega = a \xi b$. One easily finds $d\omega$:
\begin{equation} d \omega = \xi \omega + \omega \xi + a ( \tilde{\Xi} -
\xi\xi) b, \end{equation} 
and, in a special case, when  $\omega$ is of
the form $da = \xi\, a - a\, \xi$, we get: 
\begin{equation} d^2 a =
\tilde{\Xi} a - a \tilde{\Xi}, \end{equation} 
where $\tilde{\Xi}$
denotes an element in $\Om{2}$ which is the image of $\Xi$. Now, we
clearly see that division by $\CJ_2$ guarantees that $d$ is
well-defined ($d0 =0$) and $d^2 \equiv0$ if and only if we quotient out
$\CJ_3$. 

Note, that in order to obtain a differential algebra alone we do not
need to include $\CJ_1$ as the quotient.  However, we proceed so in
order to be in agreement with the {\em standard} definition of the
differential algebra for spectral triples.\Box

Before we start considering an interesting simplification, let us
say a bit more about the interpretation of $\xi\xi$ and $\Xi$.
The image of the latter (under $\pi$) is an operator which links
pairs of subspaces of $\CH$, which have the same first index $i$.
Clearly, it means that the Dirac operator ``connects'' 
$H_{ik} \to H_{pk} \to H_{ik}$, which is always the case, as 
$\xi$ is selfadjoint (for every map $H_{ik} \to H_{pk}$ 
there exists a conjugate map $H_{pk} \to H_{ik}$). Now,
$\pi(\Xi-\xi\xi)$ is a set of maps between subspaces of $\CH$
with different first index, so it does not vanish if and
only if $\xi$ connects $H_{ij} \to H_{kj} \to H_{pj}$ for
some $p \not= i$ and the composition of these maps does
not vanish.  

An interesting situation occurs when 
$\pi(\xi \tsa \xi) = \pi(\Xi)$, so that $\CJ_2 \subset \CJ_1$.

\subsubsection{Inner calculi for spectral triples.}

\begin{lem}
When $\pi(\xi \tsa \xi)=\pi(\Xi)$ then $\tilde{\Xi} = \xi\xi$ 
and the calculus remains inner in $\Om{2}$:
\begin{equation}
d\omega = \xi \omega + \omega \xi,
\end{equation}
for any one-form $\omega$. The subbimodule $\CJ$ contains the kernel of
$\pi$ and all commutators $[a,\xi \tsa \xi]$.  \end{lem}

This happens, for instance, if $\pi(\xi)\pi(\xi)$ is in $\pi(\CA)$.
In addition, if it commutes with $\pi(\CA)$ then $\CJ = \ker \pi$ and
the structure of $\Om{2}$ simplifies considerably.

The next lemma makes these observations more precise:

\begin{lem} If $\pi(\xi)\pi(\xi)$ commutes with $\pi(Z(\CA))$, where
$Z(\CA)$ is the center of $\CA$ then $\pi(\xi)\pi(\xi) = \pi(\Xi)$. The
converse is also true. \end{lem}

{\bf Proof:} Of course, if $\pi(\xi)\pi(\xi)$ commutes
with $\pi(Z(\CA))$ it commutes with every $\pi(P_i)$. Using
$\sum_i P_i \equiv 1$ we get the required identity. 
For the converse, let us assume that there exists $i$ such that
$\pi(P_i)$ does not commute with $\pi(\xi)\pi(\xi)$ and let
the commutator be $\rho$. Then it is easy to verify that
$\pi(\xi)\pi(\xi) - \pi(\Xi) = \rho$.\Box

It is quite interesting to characterize the space, 
in which $\pi(\xi)\pi(\xi)$ takes its values for this 
particular situation. We already know that it commutes with 
$\pi(\CA^o)$. Since additionally it has to commute with the center of
$\pi(\CA)$ the only possibility is that it is diagonal (i.e., when
restricted to $H_{ij}$ maps it to itself), and using the decomposition
into tensor products presented in section 2, we find that it is $
T_{ij} \otimes \hbox{id}$ on each of the subspaces $H_{ij}$. 

If the calculus (at least up to $\Om{2}$) is inner, we can
immediately make an interesting observation about the gauge
potentials.

\begin{obs} If $H$ is a selfadjoint one-form and $F(H)= dH + HH$
is its curvature form, then for the inner calculus, as described 
above, $F(-2\xi)=0$.
\end{obs}

This follows from a simple calculation. Of course, we have
presented here only the {\em formal} solution of the $F(H)=0$ equation
in $\Om{2}$. It does not guarantee that this solution is unique, it
could happen (for instance, if we take a tensor product of the spectral
triple in question  with a triple corresponding to a manifold) that
there are far more solutions or even that the whole subbimodule 
$\Om{2}$ for the discrete part vanishes. 

\subsubsection{Commutative discrete spectral triples.}

A lot more can be said about the commutative case $\CA = \C^n$. The
generators of the algebra can be identified with the 
projection operators $P_i$ introduced earlier. The structure 
of $\Omega^1(\CA)$ is completely determined by the elements 
$P_i dP_j,\; i \not=j$ (using $\sum_i dP_i=0$). In our representation 
we have:

$$ \pi ( P_i dP_j) h = P_i D P_j h - P_i \delta_{ij} D h ,$$
which for the $kl$-component becomes:

$$ \left( \pi( P_i dP_j) h \right)_{kl} = \delta_{ik} D_{kljp} h_{jp} -
 \delta_{ik} \delta_{ij} D_{klps} h_{ps} = \ldots, $$

which after taking into account the properties of $D$, becomes:

$$ \ldots = \delta_{ik} D_{iljl} h_{jl} - \delta_{ik} 
\delta_{ij} D_{klpl} h_{pl}. $$ 

Since we are interested in the case $i \not=j$ we obtain that for
$P_i dP_j \not=0 $  it is sufficient that there exists $l$ such that
$D_{iljl} \not= 0$. 

Let us now investigate $\Omega^2(\CA)$, which would be built from
the bimodule $\Omega^2_u(\CA)$ after we quotient out the subbimodule
generated by $\ker\pi \cup d {\ker \pi}$.

First, we show the representation action of the basis of
$\Omega^2_u(\CA)$, $e_{kij} = P_k dP_i dP_j$, $k\not=i$ and
$i \not=j$ (again, to see that this is the basis we use the identity
$\sum_i dP_i =0$).

\begin{obs}
For $i,j,k$ such that $k\not=i$ and $i\not=j$ $\pi(e_{kij})$ 
is a collection of maps $H_{jl} \mapsto H_{kl}$, each of the form: $$
\pi(e_{kij}) = D_{klil} D_{iljl},$$ \end{obs}

Next we shall show the structure of $d \ker\pi$.

Let us take a one-form $\rho$ in $\ker\pi$. Of course, it is
sufficient to take a generator: $\pi(P_i dP_j)=0$ for some $i\not=j$.

Let us calculate $\pi(dP_i dP_j)$. After little calculation  we get the
following result: \begin{equation} \pi (dP_i dP_j) h = - P_i D^2 P_j h.
\hbox{\ \ if } P_i dP_j \in \ker\pi, i\not=j. \end{equation}

Rewriting all that in terms of the basis of $\Omega^2_u$ we obtain that
for each $i\not=j$ if $\pi( P_i dP_j) =0$ we  have the following
generator added to the {\em junk}: 

$$ \sum_{i \not= s \not=j} P_i dP_s dP_j $$

An easy way to visualize the construction of the differential structure
coming from the commutative spectral triple is by way of graphs.
Imagine we have $n$ vertices, and each vertex is {\em split} into $n$
points. The nonzero entries of the Dirac operator shall correspond to
arrows linking two points of two different vertices or within the same
vertex (of course, due to the reality structure there exists a link
between the two).

If we look at the differential forms, the {\em fine} structure of the
vertices is less important. 

Now, as for the one forms, we know that $P_i dP_j$, $i\not=j$ exists if
there is at least one arrow connecting vertex $j$ with $i$. For a
two-form $P_k dP_i dP_j$ there must exist links $j \to i \to k$.
However, if there is no direct link $j \to k$ (note that due to chirality
this is always the case) then the sum of all forms 
$\sum P_k dP_i dP_j$ over $i$, $i \not=k$ and $i\not=j$ vanishes. This
bears a resemblance to some general structure of differential geometry on
graphs \cite{ASgr}.

\subsubsection{The metric}

On the algebra $\CA$ we have a natural trace coming from the
representation on the Hilbert space $\CH$. Since our representation
extends naturally to the differential algebra, such a trace gives a
{\em generic} scalar product on the space of one-forms. In particular,
for the generators of the commutative algebra (which is the easiest
example) we have:

\begin{eqnarray*}
\left( P_i dP_j, P_k dP_l \right) & = &\hbox{Tr} 
\left( \pi^*(P_i dP_j) \pi( P_k dP_l) \right)  \\
& = & \delta_{ik} \delta_{jl} \sum_p \hbox{Tr} (D_{lpkp} D_{kplp}).
\end{eqnarray*}

The last term of the expression is a selfadjoint map from $H_{lp}$ onto
itself, therefore the scalar product is positive. All generating
one-forms $P_i dP_j$ are orthogonal to each other and the norm of each
is fixed by the elements of the Dirac operator. Of course, should the
norm be zero then it is obvious that $\forall p \; D_{ipjp}=0$ and
hence the one-form $P_i dP_j$ vanishes.

The above choice is not unique, and one may consider its deformation of
the following type: $$ (\omega, \rho )_Z = \hbox{Tr} \left( Z
\pi(\omega) \pi(\rho)^\dagger \right), $$ where $Z$ is an operator, for
which it is sufficient to assume that $Z$ is selfadjoint, and, to guarantee
the gauge invariance, $\pi(U)Z\pi(U)^\dagger=Z$ for every unitary
$U$ in the algebra.

The reason for introducing such a scaling may follow from the requirement
that the trace has some additional symmetry (for instance we might fix
it to recover the Haar measure on the algebra, when it has a Hopf
algebra structure).

Of course, changing the trace would have significant physical 
consequences.  For details of the construction and application to the
Standard Model, see the forthcoming paper \cite{MP}. 

\subsection{Finite spectral triples with real algebras.}

So far we have discussed finite spectral triples with algebras over
$\C$. Since the possible application of noncommutative geometry
in high energy physics  uses rather a real algebra, we shall now
briefly discuss the above construction assuming that $\CA$ is over the
field of real numbers.

\begin{lem}[See \cite{EMS}]
The simple finite algebras over $\R$ are of the type $M_n(F)$, where
$F$ can be $\R,\C$ or $\H$. \end{lem}

Therefore every semi-simple algebra over $\R$ can be decomposed as:
\begin{equation} \CA = \bigoplus_i M_{n_i}(F_i), \end{equation} where
each $F_i$ could be the field of real, complex or quaternionic numbers.

We shall now investigate irreducible representations of $\CA$ on
a {\em complex} Hilbert space. If one considers $M_n(\R)$ and $M_n(\H)$ 
as real algebras, then each of these has only one irreducible representation on $\C$
(for quaternions the conjugate representation is equivalent to the
fundamental one).  In the case of  $M_n(\C)$ we have two inequivalent irreducible
representations: $$ \pi(a) = a \;\;\;\; \bar{\pi}(a) = \bar{a} $$ where
$\bar{a}$ means the involution of matrix elements. 

How does the description of our spectral triple change? Apparently for
$M_n(\R)$ and $M_n(\H)$ there is no change at all, whereas for
$M_n(\C)$ we have to take into account the existence of these two
inequivalent representations: \begin{lem} If $F_i=\C$ then the Hilbert
subspace $H_{ij}$ decomposes into $H_{ij}$ and $H_{\bar{i}j}$, with the
representation of $M_{n_i}(F_i)$ by $m$ and $\bar{m}$, respectively.
\end{lem}

Of course, an analogous decomposition takes place if one takes into
account the right action of the algebra on $\CH$.

Until now, there has not been much change apart from the additional
splitting of the Hilbert space. All other results remain unchanged, in
particular, $\gamma$ cannot be different for $H_{ij}$ and
$H_{\bar{i}j}$.

If we look at the restrictions for the Dirac operator, which come
from (\ref{D-prop}) we see that the constraints on allowed values
of $D$ also remain unchanged, however, one should add some more,
which arise from the fact that one cannot mix the fundamental 
and the conjugated representation. To see this explicitly,
in the proof of Lemma 6 one should replace , $b=P_j$ by $b=z P_j$ 
for any complex number~$z$.  

Generally speaking, the derivation of the constraint which we had for
the complex case, could be repeated separately for each of the real
situations. Therefore, we may summarize the result by stating that $D$
does not mix two inequivalent representations of the same algebra.

\subsection{The intersection form.}

The intersection form in $K$-theory is a map $K_*(\CA) \times
K_*(\CA) \to \Z$ and its invertibility (Poincar\'e duality)
is fundamental for a characterization of homotopy types of spaces,
which possess the structure of a smooth manifold.  

For finite complex spectral triples the intersection form, evaluated on
generators of the $K_0$-group, $e, f$, is \cite{Connes}:

$$ \langle e, f \rangle = \langle  D, e \otimes f^o \rangle, $$

and, generally, the pairing $\langle D, E \rangle$ for a projector $E$
means:

$$ \langle D, E \rangle = \hbox{dim\ }\hbox{coker}_{E\CH_R}  (ED^+E) - 
 \hbox{dim} \hbox{ker} (ED^+E) $$

where $\CH_L = \frac{1}{2}(1+\gamma) \CH$,
$\CH_R = \frac{1}{2}(1-\gamma) \CH$ and $D^+ = \frac{1}{4} 
(1-\gamma)D(1+\gamma)$.

Since for matrix algebras the group $K_1$ is trivial (any element can
be deformed to $1$), the only nontrivial part is $K_0$. For $M_n(\C)$
all projectors are equivalent (within $M_\infty(M_n(\C))$) to a
diagonal matrix with $1$ in the first diagonal entry and zeroes
elsewhere, which we shall call $e$. Of course, for a direct sum of
$N$-simple algebras the corresponding $K$-theory will have $N$
independent generators, each coming from the simple component of the
sum.

Let us calculate the $\langle e_i, e_j \rangle$ entry. 
First, we have to determine the dimensions of Hilbert subspaces
appearing in the above formulae: \begin{eqnarray} & \hbox{dim} \; e_i
e_j^0 \CH_L = \left\{ \begin{array}{ll} r_{ij} & \hbox{if\ }
\gamma_{ij} = -1 \\ 0 & \hbox{if\ } \gamma_{ij} = 1 \end{array} \right.
\\ & \hbox{dim}\;  e_i e_j^0 \CH_R = \left\{  \begin{array}{ll}  r_{ij}
& \hbox{if\ } \gamma_{ij} = 1 \\
 0 & \hbox{if\ } \gamma_{ij} = -1 \end{array} \right. 
\end{eqnarray}
Now, observe that out of $e_i e_j^o \CH_L$  and $e_i e_j^o \CH_R$ one must
always be empty, so the $e_i e_j^o D^+ e_i e_j^o $ operator, which acts
between them has either empty domain or empty target space. Thus, the index
is independent of $D$, and reads: $$ \langle e_i, e_j \rangle =
\gamma_{ij} r_{ij} = q_{ij} $$ so we recover the matrix $q$ as our
intersection form.

For real algebras the calculation is slightly more complicated. First,
one has to take into account the form of the projectors for quaternions,
which results in the doubling of the corresponding entries in the
intersection form. Moreover, for real algebras the group $K_1$ is
nontrivial, as $K_1(M_n(\R))=\Z_2$, however, this still does not change
the intersection form if one does not take into account torsion in
$K$-theory.   

\section{Hopf algebra symmetries.}

Spectral triples are expected to be an extension of the notion of
Riemannian manifolds to noncommutative geometry. Similarly as in the
standard differential geometry we may attempt to explore symmetries of
such manifolds. We expect, however, that the correct idea of a symmetry
should be extended from a group to a Hopf algebra (for a good
introduction to Hopf algebras and quantum groups we refer the reader to
\cite{Majid}). Another potentially interesting point is the speculation
\cite{Connes,CCh,Co2} that the finite algebra of the Standard Model
originates from the $q$-deformation of the Spin structure and may
possess some Hopf algebra symmetry, for details see also \cite{Coq}.

We shall try to answer the question whether spectral triples admit (and
if so, to what extent) a symmetry understood as a Hopf algebra
symmetry. Of course, the natural thing is to use the Haar measure as a
trace on the algebra and to extend it to the whole representation
space. Another possibility is to demand that differential structures of
the spectral triple are also symmetric, for instance, that they are
left-covariant or bicovariant. Furthermore, the Hilbert space may be
investigated for the existence of a comodule structure.  Last not
least, we may turn our attention to the tensor product of
representations (this is allowed only if a coproduct on the algebra
exists), which physically has the meaning of constructing  composite
states. In this paper we concentrate our efforts on the first two
aspects of symmetries, leaving the others for future investigations.

\subsection{Discrete group structure on $\C^n$ algebra}

One of the simplest examples of real spectral triples are
those based on a finite-dimensional commutative algebra $\C^n$
understood as function algebras on a discrete $n$-point space.
The general construction scheme was presented by Connes \cite{Connes}
and falls into  the general classification of complex spectral triples
presented in the previous section.

Suppose  that we assume the discrete space
to have the structure of a finite group and demand that the 
first-order differential calculus generated by $D$ is 
bicovariant. 

First, let us recall the basics of the bicovariant differential
structure for finite groups (details in \cite{AH}). The calculus is
generated by left-covariant forms: $$ \chi^g = \sum_h  e_{gh} d e_g$$
where $e_g$ is a Kronecker-delta function: $1$ at $g$ and zero
elsewhere. The right coaction on them is: $$ \Delta_R \chi^g = \sum_h
\chi^{hgh^{-1}} \otimes e_h $$ Using the spectral triple approach, we
may calculate the $\chi^g$ forms as operators on our Hilbert space. It
appears that $$ (\chi^g h)_{ij} =   D_{ij,(ig^{-1})l} h_{(ig^{-1})l}.
$$ Because of the properties of $D$ we see that it is only possible
that $j=l$ and therefore finally: $$ (\chi^g h)_{ij} =  
D_{ij,(ig^{-1})j} h_{(ig^{-1})j}, $$ so it is a collection of operators
from $H_{(ig^{-1})j} \to H_{ij}$. Of course, if $\gamma_{ij}
\gamma_{(ig^{-1})j} = 1$ they all vanish because $D$ acts only between
spaces of different chirality.

What we have obtained so far is the representation of left-invariant
forms as operators on our Hilbert space. We see that all $\chi^g$,
which do not vanish are linearly independent (since they act between
different Hilbert subspaces). 

If we demand that the calculus is bicovariant we immediately see that
if for any $g$, $\chi^g$ vanishes, then for every $h$ $\chi^{hgh^{-1}}$
must vanish as well. Of course, for commutative groups such a condition
is void, however, for nonabelian ones it leads to certain restrictions
on possible forms of $D$, which we make precise in the following lemma:

\begin{lem} If the spectral-triple calculus on $\C^n$ with a structure
of the discrete group $G$ is bicovariant then: \begin{eqnarray*} &
\hbox{if\ }\chi^g \not=0 \hbox{\quad \em then\ } \forall i \in G \;
\exists j \in G: D_{ij,(ig)j} \not= 0 \\ & \hbox{if\ }\chi^g =0
\hbox{\quad \em then\ } \forall i \in G \; \forall j \in G:
D_{ij,(ig)j} = 0 \end{eqnarray*} \end{lem}

The proof is an immediate consequence of previous observations.

\subsubsection{Function algebra on $S_3$.}

Consider as an example the $6$-dimensional space $S_3$. We shall check
what restrictions on the spectral triple and the Dirac operator are set
by the requirement that the corresponding calculus be bicovariant.

First, let us take the universal calculus. It appears that there
exists a spectral triple which gives such a calculus.
The simplest example is given by the following intersection form:
$$
\left( \begin{array}{cccccc}
-1 & 1 & 1 & 1 & 1 & 1 \\
1 & -1 & 1 & 1 & 1 & 1 \\
1 &  1 &-1 & 1 & 1 & 1 \\
1 &  1 & 1 &-1 & 1 & 1 \\
1 &  1 & 1 & 1 &-1 & 1 \\
1 &  1 & 1 & 1 & 1 &-1 \end{array} \right)
$$
and is $36$-dimensional! It is quite easy to see that this must be the
case - for each $\chi^g$ not to vanish, the Dirac operator must link
(for every $i$) some $H_{(ig^{-1})j}$ with $H_{ij}$. For the
intersection matrix it means that for every row there always exists an
entry which has an opposite sign to an entry in the same column of
another row. It is easy to verify that such a matrix is nondegenerate so
that the Poincar\'e duality axiom is fulfilled. Of course, to get the
universal calculus, all possible Dirac operator elements must not
vanish.

A more interesting case is when we require that the calculus is the
lowest-dimensional bicovariant. For $S_3$ such a calculus is generated
by $\chi^{ab}$ and $\chi^{ba}$, where $a,b$ denote the generators of
$S_3$. 

If we order the elements of $S_3$ as $a,b,c,1,ab,ba$, the
example of a spectral triple which gives the desired 
calculus is defined by:
$$
\left(
\begin{array}{cccccc}
1  & 1 & -1 & 0 & 0 & 0 \\
1  & -1 & 0 & 0 & 0 & 0 \\
-1 & 0 & 0 & 0 & 0 & 0 \\
0  & 0 & 0 & 1 & 1 & -1 \\
0  & 0 & 0 & 1 & -1 & 0 \\
0  & 0 & 0 & -1 & 0 & 0 \end{array} 
\right).
$$

Again, all the possible elements of the Dirac operator (components
which act between spaces of different chirality) must not vanish. Of
course if we slightly change the intersection form (introduce, for
instance, other nonzero entries in the intersection matrix) the
bicovariance (2-dimensional) shall also be preserved provided that the
additional allowed Dirac operator components vanish!

For physical interpretation this would mean that the symmetry requires
some of the fermion particles in the model to be massless, whereas in
the spectral triple presented above all mass matrices must be nonzero.

Noncommutative finite algebras cannot carry such a simple group
structure and therefore we have to look for nonabelian Hopf
algebras. This in itself is an interesting topic and very little is known
about the general classification of such (semisimple) objects. The
easiest examples come from group algebras and as we try to make our
examples  comprehensible, we shall use the simplest one of them, the
group algebra $\C{}S_3$. 

\subsection{Group algebra $S_3$}

We shall begin here with some generalities on Hopf and bicovariant
differential structures on group algebras and then we shall present the
interesting example of $S_3$.

\subsubsection{Generalities on group algebras.}

The group algebra $\C{}G$ has the natural structure of a Hopf algebra
with the following coproduct, counit and antipode: $$ \Delta g = g
\otimes g, \;\;\; \epsilon(g) =1, \;\;\; Sg = g^{-1} $$

The adjoint coaction is trivial:
$$ ad(g) = g \otimes 1 $$

The Haar invariant measure $\mu$ is just $\mu(g)= 0$ for $g$ different
from the neutral element of the group.

Finally let us state some observations on the differential calculi.

\begin{obs} 
The differential calculi on group algebras are always inner, i.e.,
there exists a one-form $\chi$, such that $dg = [g, \chi]$ \end{obs}

It is sufficient to take:
$$\chi = - \frac{1}{N_G} \left( \sum_h h^{-1} dh \right). $$
Then:
\begin{eqnarray*}
&  [g, \chi ]= - \frac{1}{N_G} \sum_h \left( g h^{-1}\, dh - h^{-1}
dh\, g \right) \\ & = - \frac{1}{N_G} \sum_h \left( gh^{-1}\, dh -
h^{-1} d(hg) - dg \right) \\ & = \frac{1}{N_G} \sum_h dg = dg
\end{eqnarray*}

\begin{obs}
Since $\C{}G$ is cocommutative every left-covariant calculus is
bicovariant. There exists a 1:1 correspondence between bicovariant
calculi and representations of the group $G$ \end{obs}

The first remark is obvious, so let us concentrate on the latter.
Suppose we have a representation of $G$ on some vector space $V$. We
may define $\Omega_1$ as a right-module $\C{}G \otimes V$, with the
multiplication from the right by elements of $\CA$ only on the first
component of the tensor product. Then the following rule (it is
sufficient to define it only for the generators) makes $\Omega_1$ a
bimodule:

$$ h (g \otimes v)   = (hg) \otimes (\pi(h) v). $$

Of course, $\Omega_1$ is a bicovariant bimodule, with all the elements
of type $1 \otimes v$ being left and right invariant. We already know
that every differential calculus is inner and, moreover, it is easy to
check that the one-form $\chi$ is left and right invariant Therefore,
to construct the calculus, it is sufficient to choose a $\chi = 1
\otimes v \in 1 \otimes V$ and define $dg = [g, \chi]$. It remains to
show that $\hbox{Im}\, d$ generates $\Omega_1$. This is equivalent  to
verifying  whether $\{ v - \pi(h) v \}, h \in G$ generates the whole
space $V$. Note that the representation $\pi$ could be linearly
extended on the whole group algebra $\C{}G$. Then the latter
requirement tells us that the action of the kernel of the counit 
on $v$ spans the whole space $V$. 

The proof of the converse (every bicovariant calculus provides us with
a representation of the group) is obvious.

\begin{obs} The central part of $\Om{1}$ consists of all one-forms of
the type: $ \omega = \sum_g g \otimes v_g $, where $v_g \in V$ are
vectors satisfying for every $g,h \in G$: \begin{equation} \pi(h) v_g =
v_{h g h^{-1}}. \label{eq34} \end{equation} \end{obs}

Of course, this condition may be void for some representation of
the group $G$ (in particular for a nontrivial one-dimensional
representation).

{\bf Proof:} Every one form can be written as $\omega=\sum_g g \otimes
v_g$ for some $v_g \in V$. If we require that $h\omega=\omega h$ for
every $h \in G$ then it is easy to verify that (\ref{eq34})
follows.\Box

As an example we shall consider the smallest noncommutative group
algebra $\C{}S_3$.

\subsubsection{The group algebra $\C{}S_3$}

The group $S_3$ has two generators $a,b$ with the following
multiplication rules:

$$ a^2 = e, b^2=e, aba=bab$$ 

It is convenient to name the element $aba=c$, $c^2=e$ and use the
rules for multiplication between $a,b,c$:

$$ ab = bc, \;\;\;\; ac = ba, \;\;\;\; bc = ca  $$

As an algebra $\C{}S_3$ is isomorphic to the algebra 
$M_2(C) \oplus \C \oplus \C$. We shall present the form of this
isomorphism $i$ for the generators $a,b$:

$$ i(a) = \left( \begin{array}{cccc} 1 & 0 & & \\
0 & -1 & & \\ & & 1 & \\ & & & -1 \end{array} \right), \;\;\;\;
 i(b) = \left( \begin{array}{cccc} -\frac{1}{2} & \frac{\sqrt{3}}{2} &
 & \\
\frac{\sqrt{3}}{2}  & \frac{1}{2} & & \\ & & 1 & \\ & & & -1
\end{array} \right), $$

The Haar measure on the Hopf algebra is defined as a linear map, which
satisfies:

$$  1 \mu(f) = (\hbox{id} \otimes \mu) \Delta f = 
(\mu \otimes \hbox{id}) \Delta f. $$

On any group algebra there exists only one such measure, which is
$ \mu(g) = 0, g \not=0$ and $\mu(1)=1$ (normalized). For our 
$4$-dimensional representation ($M \in M_2(C), p,q \in C$) 
it is given by:

$$ \mu \left( \begin{array}{ccc} \hbox{\huge M} & & \\
 & p & \\ & & q \end{array} \right) = \frac{1}{3} \hbox{Tr\ } M +
\frac{1}{6} p + \frac{1}{6} q $$

Before we start constructing spectral triples let us have a look
at the covariant differential calculi, constructed according to
the general scheme. Therefore, first we have to give the
representations of the group $S_3$:

\begin{center}
\begin{tabular}{|c|c|c|}
\hline
dimension & $\pi(a)$ & $\pi(b)$ \\ \hline
1* & $1$  & $1 $ \\ \hline
1 & $-1$ & $-1$ \\ \hline
2 & $\displaystyle \left( \begin{array}{cc} 1 & 0 \\ 0 & -1 \end{array}
\right)$
  & $\displaystyle \left( \begin{array}{cc} -\frac{1}{2} &
  \frac{\sqrt{3}}{2}  \\
\frac{\sqrt{3}}{2}  & \frac{1}{2}  \end{array} \right)$ \\ \hline
\end{tabular}
\end{center}
where 1* is the trivial representation. 
Every representation of higher dimension is a direct sum of 
them.

\begin{itemize}
\item  {\bf 1-dimensional calculus}
For the one-dimensional calculus we take the representation $1$ 
(the trivial one gives $d\equiv0$).

There is only one generating form $\chi$ and the bimodule rules
follow from the one-dimensional representation of $S_3$. Since
there is only one such representation, $a=b=-1$ we get:
$$a \chi = -\chi a, \;\;\;\;\; b \chi = - \chi b.$$ 
and the external derivative is
$$ d a = a \chi, \;\;\;\;\; db = b \chi,  \;\;\;\;\; dc = c \chi, $$
notice that: $$ d(ab) = 0, \;\;\;\;\; d (ba) = 0. $$ Such a calculus is
extendible to a higher-order calculus with $d^2=0$ provided that we set
$d \chi = -\chi \chi$.

This calculus has an interesting splitting property:
\begin{obs}
Let $A = M_2(C)$ and $B=\C \oplus \C$, $\C{}S_3 = A \oplus B$.
Then the split short exact sequence of algebras:
$$ 0 \to A \to \C{}S_3 \to B \to 0 $$
extends to a split short exact sequence of differential modules:
$$ 0 \to \Omega(A) \to \Omega(\C{}S_3) \to \Omega(B) \to 0 $$
where $\Omega(\C{}S_3)$ is the one-dimensional bicovariant calculus,
and $\Omega(A)$, $\Omega(B)$ are its restriction to the corresponding
subalgebras. \end{obs}

This tells us immediately that no spectral triple data could generate
this calculus. Indeed, since the Dirac operator $D$ can act only
between spaces of different left (or right) representation of the
algebra, such a splitting of differential bimodules is not possible.

\item {\bf two-dimensional calculus}
As we have pointed out there exists only one nontrivial two-dimensional
representation of $S_3$ (here nontrivial means that that there exists a
vector $v$ such that $(\ker \epsilon) v$ spans the whole representation
space. For instance,  this does not hold for $1 \oplus 1$). Various
possibilities for the choice of $v$ would correspond to differential
calculi related by automorphisms of the group algebra. 

If we choose the above matrix form of $a$ and $b$ and
$v= (1,-\sqrt{3})$ we arrive at the following relations: 
$$ \chi^1 a = - a \chi^1, \;\;\;\; \chi^2 b = - b \chi^2, $$
$$ \chi^1 b = b (\chi^1-\chi^2), \;\;\;\; \chi^2 a = a( \chi^2 -
\chi^1), $$ where $\chi^1= a\, da$ and $\chi^2 =b \, db$.

The same relations expressed in terms of $dg$ become
$$ (da) b = c (da) - a (db) \;\;\;\; (db) a = c (db) - b (da) $$
then $ dc = (da) ba + a (da) b + ab (da) = 0$.

All other two-dimensional differential calculi could 
be obtained through an automorphism of the algebra.

Note that for this calculus we may construct an one form in the
center of $\Om{1}$. For instance, the vector 
$a v_a + b v_b - c ( v_a +v_b)$, where $v_a, v_b$ are eigenvectors (to
eigenvalue $1$) of $\pi(a),\pi(b)$, respectively, gives such a form.
Using $\chi^a$ and $\chi^b$ we can write it as $$ (2a-b-c)\chi^2 +
(2b-a-c)\chi^1,$$

so, again this rules out the compatibility between bicovariance and
spectral-triple calculus for $\C{}S_3$. 

Our result extends also for higher dimensional calculi, as all
higher order representations of $S_3$ contain the 
two-dimensional part (and, if not, they must be composed out of 1 and
1*) and therefore one may repeat the arguments. \end{itemize}

Thus, we may state a simple no-go lemma:

\begin{lem}
For $\C{}S_3$ the differential structure from the spectral triple
construction cannot carry a bicovariance symmetry.
\end{lem}

Knowing that we may now concentrate on spectral triples and
the symmetry of the measure.

\subsection{Spectral Triples for $\C{}S_3$}

Following the general construction scheme presented earlier we shall
now turn our attention to an example of low-dimensional spectral
triples for the algebra $\C{}S_3$. We will be interested in the
non-trivial case, which admits a non-zero Dirac operator.

It is easy to verify that the lowest dimensional spectral triple is
given by the intersection form:

$$ \left( \mth{0}{1}{0}{1}{0}{-1}{0}{-1}{1} \right). $$

The Hilbert space has dimension $7$ (it is odd, as it 
contains a vector in the subspace $H_{33}$ such that 
$Jv=v$!).

Instead of the Dirac operator it is more convenient to use
the one-form $\xi$. Since the Dirac operator $D$ can only have 
nonzero entries $D_{12,32}$ and $D_{23,33}$ (and, respectively, their
hermitian conjugates and $J$-conjugates) we have the following $\xi$
(where the Hilbert spaces on which $\xi$ acts are ordered as follows:
$H_{12}, H_{32}, H_{23}, H_{33}$): $$ \pi(\xi) = \left(
\begin{array}{ccccc} 0 & 0 & x  & 0 & 0\\ 0 & 0 & y  & 0 & 0 \\ x^* &
y^* & 0 & 0 & 0 \\ 0 & 0 & 0  & 0 & z \\ 0 & 0 & 0  & z^* & 0
\end{array} \right) $$

\begin{obs} 
The operator $\pi(\xi)\pi(\xi)$ is in $\pi(\CA)$
provided that $xx^*+ yy^* = zz^*$. In any case it commutes
with the image of the center of the algebra and 
therefore it follows from Lemma 11 that $\pi(\xi)^2=
\pi(\Xi)$ and the calculus is inner: 
$$ \pi(\xi)^2 = 
\left( \begin{array}{ccccc}
xx^* & xy^* & 0  & 0 & 0\\
yx^* & yy^* & 0  & 0 & 0 \\
0  & 0 &  xx^*+yy^* & 0 & 0 \\
0 & 0 & 0  & zz^* & 0 \\
0 & 0 & 0  & 0 & zz^* \end{array} \right).$$
\end{obs}

We shall finish the investigation by writing the {\em measure
operator}, which gives the Haar measure on $\C{}S_3$.

\begin{obs}
If $Z = \frac{1}{3} P_1 + \frac{1}{18} P_2 + \frac{1}{12} P_3$, 
where $P_i$ denote the corresponding projections in the 
algebra, then for every $a \in \CA$ 
$$\hbox{Tr}_Z \pi(a) = Tr\left( \pi(Z)\pi(a) \right), $$
is the Haar measure on $\C{}S_3$.
\end{obs}

One may now use this deformed scalar product on forms to
calculate the Yang-Mills action. Note that
$[D,\pi(Z)] \not= 0$ and the results would differ from the ones
obtained using the {\em generic} scalar product. For physical
models this would result in different  mass relations between 
the gauge bosons in the theory based upon such a triple.  

\section{Conclusions}

The world of noncommutative geometry is much bigger than that
of classical geometry.  In classical differential geometry, once
given a smooth manifold we already know its differential bundle. Given
a Lie group we have a natural notion of the Lie group action on
differential structures. All this seems to be lost when we enter the
noncommutative world. For a single algebra we have many choices of
differential structures and  there is no universal rule to tell us
which one to choose. Symmetries, or at least  group symmetries are in
many cases no longer present or drastically reduced.

The notion of spectral triples is an attempt to bring some order
into our understanding of noncommutative differential geometry,
by specifying what structures must appear in the models to make
them a real {\em geometry}.

What we investigate in this paper are the restrictions,
that this order introduces to the realm of $0$-dimensional
geometry.  Even more, we attempt to answer the question whether the
proposed structures admit and restrict the possible $0$-dimensional
noncommutative symmetries.

We classify finite spectral triples and the allowed Dirac
operators. This, for instance, enables us to state that the Standard
Model Dirac operator is not the most general one. In fact, one may
consider an identical spectral triple with some additional components
of the Dirac operator, which link left-handed antileptons with
right-handed quarks. Of course, the assumption that such a component
exists might have some deep consequences for the model, as it could
potentially lead to the breaking of the $SU(3)$ strong symmetry.
Another consequence (mentioned originally by D.Testard, see
\cite{Test}) is that Poincar\'e duality enforces the absence
of right-handed neutrinos. This is correct provided that we add a pair of
them (particle and antiparticle). However, it is easy to verify that
adding only one right-handed Majorana particle (we doubt whether the
name {\em neutrino} would be justified) still does preserve the
Poincar\'e duality. These and other observations concerning the
Standard model will be the topic of our forthcoming paper \cite{PS2}.

The symmetries of spectral triples are a far more complicated problem.
We have seen that only in some cases there exists an extension (on the
differential level at least) of the Hopf algebra structure into the
spectral triple theory. Whether this is a general pattern is difficult
 to say, however,  we believe that symmetries can also be realized in
 the 
$0$-dimensional 
case.  The examples we have analyzed do not give a conclusive answer,
what they suggest, however, is that the symmetry restrictions could be 
much stronger and of different type than in classical geometry.

\ \\

{\bf Acknowledgments.} We would like thank Florian Scheck for
valuable discussions and J.M.Warzecha for remarks on the manuscript. \\
\ \\ While preparing a final version of this paper we learned that
similar results were obtained by T.Krajewski \cite{Kr}, whom we would
also like to thank for comments and remarks.

\end{document}